\documentclass[conference,a4paper]{IEEEtran}
\usepackage{xcolor}
\usepackage{balance}

\usepackage[american]{babel}

\usepackage{cite}
\usepackage{xcolor}
\usepackage{balance}

\ifCLASSINFOpdf
  \usepackage[pdftex]{graphicx}
\else
  \usepackage[dvips]{graphicx}
\fi
\usepackage{calc,ifthen,intcalc}
\usepackage{tikz}
	\usetikzlibrary{arrows,shapes,positioning,calc,fadings,patterns,decorations}
	\usetikzlibrary{decorations.pathreplacing,decorations.pathmorphing}			
	\usetikzlibrary{arrows,arrows.meta}
	\usetikzlibrary{shapes.geometric}

\usepackage{amsmath}
\usepackage{bm,amssymb}
\ifCLASSOPTIONcompsoc
 \usepackage[caption=false,font=normalsize,labelfont=sf,textfont=sf]{subfig}
\else
 \usepackage[caption=false,font=footnotesize]{subfig}
\fi

\usepackage{stfloats}
\newcommand{\latexConstant}{28.45}
\newcommand{\xOffset}{0.9}
\newcommand{\xComplement}{0.1}
\newcommand{\xScale}{0.4}
\newcommand{\yScale}{\xScale*\xOffset}

\usepackage{lipsum}
\usepackage{booktabs}

\newcommand{\bs}[1]{\bm{\mathrm{#1}}}
\newcommand{\abs}[1]{\left|#1\vphantom{f}\right|}

\renewcommand{\d}{\ensuremath{\partial}}

\newcommand\figref[1]{Fig.~\ref{#1}}

\graphicspath{{images/}}

\usepackage{lipsum}

\begin{document}

%\title{Verification and Simplified Approximation of Optimized Reconfigurable Intelligent Surfaces}

%\title{Convex Optimization of Reconfigurable Intelligent Surfaces: Verification by Measurements}
\title{Validating Convex Optimization of Reconfigurable Intelligent Surfaces via Measurements}

\author{\IEEEauthorblockN{
Hans-Dieter Lang\IEEEauthorrefmark{1},   % 1st author, 1st affiliations
Michel A. Nyffenegger\IEEEauthorrefmark{1},   % 2nd author, 2nd affiliations
Sven Keller\IEEEauthorrefmark{1},    % 3rd author, 3rd affiliations
Patrik Stöckli\IEEEauthorrefmark{1},\\    % 3rd author, 3rd affiliations
Nathan A. Hoffman\IEEEauthorrefmark{1},
Heinz Mathis\IEEEauthorrefmark{1},    % 3rd author, 3rd affiliations
and Xingqi Zhang\IEEEauthorrefmark{2}      % 4th author, 4th affiliations
}                                     % ...
%\\
\IEEEauthorblockA{\IEEEauthorrefmark{1}% 1st affiliations
OST – Eastern Switzerland University of Applied Sciences, Rapperswil, SG, Switzerland. Email: hansdieter.lang@ost.ch}
\IEEEauthorblockA{\IEEEauthorrefmark{2}% 2nd affiliations
School of Electrical and Electronic Engineering, University College Dublin, Dublin, Ireland.} 
% \IEEEauthorblockA{ \emph{\{hansdieter.lang,michel.nyffeengger\}@ost.ch} }
}



\maketitle

\begin{abstract}
Reconfigurable Intelligent Surfaces (RISs) can be designed in various ways. A previously proposed semi\-definite relaxation-based optimization method for maximizing power transfer efficiency showed promise, but earlier results were only theoretical. This paper evaluates a small RIS at 3.55\,GHz, the center of the 5G band “n78”, for practical verification of this method. The presented results not only empirically confirm the desired performance of the optimized RIS, but also affirm the optimality of the resulting reactance values. Additionally, this paper discusses several practical aspects of RIS design and measurement, such as the operation of varactor diodes and time gating to omit the direct line-of-sight (LOS) path.%\\[-3.5mm]
\end{abstract}

\vskip0.5\baselineskip
\begin{IEEEkeywords}
 Reconfigurable Intelligent Surface (RIS), Metasurface, Scattering,   Optimization Methods, Loaded Antennas.
\end{IEEEkeywords}

% \fi
% \IEEEpeerreviewmaketitle

\vspace*{4pt}
\section{Introduction}
% no \IEEEPARstart“n78”

\subsection{Background}

%Reconfigurable intelligent surfaces (RISs) have been drawing attention across various research domains including applied electromagnetics, antennas \& propagation, and mobile communications~\cite{hu18,huang19,renzo20}. They are broadly considered to be pivotal in upcoming wireless communication technologies and networks, including those categorized as "beyond 5G" (B5G) and extending into 6G and beyond \cite{huang20,saad20,liu21}.

Reconfigurable Intelligent Surfaces (RISs) are broadly considered to be pivotal in upcoming wireless communication technologies and networks, including those categorized as “beyond 5G” (B5G) as well as 6G and beyond \cite{huang20,saad20,liu21}. Hence, they have been drawing attention across various research domains including applied electromagnetics, antennas \& propagation, and mobile communications~\cite{hu18,huang19,renzo20}. 

The underlying principle is to redirect wireless signals, which may otherwise be lost or “spilled”, back into the system, thereby enhancing both coverage and throughput of wireless connections. To ensure these reflections prove advantageous, reflectors must be both adaptable and “intelligent”. This requires them to be controllable in order to optimize overall performance for each respective link. Moreover, such optimization should consider all users and the entire channel and incorporate all additional scatterers in the vicinity. 

%This should be done considering all users and the full channel, including all additional scatterers in the vicinity. 
%Such intelligent reflectors aim to diminish the necessity for extra base stations, thereby conserving costs and energy. 
%Employing RISs should improve system performance relative to cost and boost energy efficiency and improve sustainability overall in future wireless communication systems. They aim to diminish the necessity for extra base stations, thereby not only conserving costs and energy, but also contributing to a more environmentally friendly and sustainable network infrastructure by minimizing excessive use of materials.
%Consequently, employing RISs should improve system performance relative to cost and also boost energy efficiency and sustainability overall in future wireless communication systems. 

Ultimately, the employment of RISs strives not only to mitigate the necessity for supplementary base stations, conserving both costs and energy, but also to amplify energy efficiency and fortify overarching sustainability in forthcoming wireless communication systems.

\subsection{Overview of Related Work}

Several RIS design methodologies are explored in the literature, see e.g., \cite{huang19,renzo20}.
Binary control of RIS elements is a common method that simply switches elements on or off, providing a simple but limited control mechanism. 
Other prevalent approaches incorporate impedance boundary considerations or involve typical metasurface design methods, frequently relying on (locally) periodic approximations that also have their limitations~\cite{li23}.

This project adopts reactive loading as in the design of reconfigurable reflectarrays~\cite{hum14}. This enables continuous phase adjustment of the reflected wave, offering enhanced RIS control, without the limitations of binary on/off states. The question is, how the reactances are obtained.

In \cite{lang23} it was proposed to use an optimization method to maximize power transfer from a transmitter via RIS to a receiver to find the optimum loading reactances. The method is based on a semidefinite relaxation of the otherwise non-convex underlying optimization problem that showed promising results in simulation, however an empirical verification via measurements was amiss. It had only been verified in relatively simple, non-reconfigurable transmit array scenarios~\cite{nyffenegger21}. Optimality of the loading reactances was not shown. This paper aims to fill that gap, thereby confirming the validity and practicality of the proposed method.

\subsection{Objectives \& Outline}
This paper delivers empirical validation via measurements of an optimization framework previously proposed for RIS, which to this point had only been explored in simulation. To facilitate this, an RIS design is described and a comprehensive measurement setup is developed. Several measurement results will affirm the validity and applicability of the aforementioned optimization method. Furthermore, various practical facets of both the RIS design and the measurement approaches are explored and discussed.

In the following, Section II delves into the specific wireless scenario under consideration, while Section III explores the practical implementation of the reconfigurable reflectarray, accompanied by a discussion of its inherent limitations. Finally, Section IV presents measurement results, serving to validate the proposed approach.

\vspace*{4pt}
\section{Preliminaries}%Application Scenario and Review of Optimization Method}

\subsection{Wireless Channel Setup}
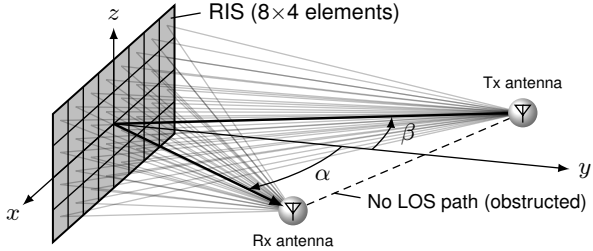
\begin{figure}[!htb]\centering
	%\vspace*{-1mm}
	\begin{tikzpicture}[scale=1,>=latex,line width=0.5pt]\sffamily\small
		%\node[opacity=0.66] at(0,0){\includegraphics[width=4cm]{principle.jpg}};
		
		\newcommand{\shiftX}{0*\latexConstant}
%%%%% F R O N T %%%%%
\begin{scope}[yshift=\shiftX*-\xComplement,xshift=\shiftX]
\begin{scope}[every node/.append style={yslant=\xOffset,xscale= \xScale},yslant=\xOffset,xscale= \xScale,xshift=0,yscale=0.85]
	
	\draw[line width=0.75pt,fill=lightgray] (-2,-1)rectangle(2,1);
	\draw (-2,0.5)--(2,0.5);
	\draw (-2,0)--(2,0);
	\draw (-2,-0.5)--(2,-0.5);
	\foreach \x in {-1.5,-1,...,1.5}{
	\draw (\x,-1)--(\x,1);
	}
	
	\sffamily
	
	\draw[->] (0,0)--(-3,0)coordinate[below left=5pt](xtemp);
	\draw[->] (0,0)--(0,1.5)coordinate[above=5pt](ztemp);

\end{scope}
\end{scope}

\renewcommand{\shiftX}{0*\latexConstant}

%%%%% S I D E %%%%%
\begin{scope}[xshift=-\xScale*\shiftX,yshift=-\yScale*\shiftX]
\begin{scope}[every node/.append style={yslant=\xOffset-1},yslant=\xOffset-1]	

	%\draw[blue!90!black,densely dashed] (-0.215,-4)--(-0.215,-0.6) (0.213,-3.98)--(0.213,-0.6);
	
		\draw[->] (0,0)--(6,0)node[right]{$y$};
	\node at(xtemp){$x$};
	\node at(ztemp){$z$};

\end{scope}
\end{scope}

\def\aalpha{35}
\def\bbeta{-20}

%%%%% T O P %%%%%
\begin{scope}[yshift=0*\latexConstant,scale=-2]
\begin{scope}[every node/.append style={yslant=\xOffset,xslant=-1,xscale= \xScale},yslant=\xOffset,xslant=-1,xscale= \xScale]

\def\delta{0.1758}
		\def\T{2.5}
		\def\R{2}
		
		%\draw[->] (-0.65,0)--(0.65,0)node[right]{$x$};
		
		\draw[line width=0.9pt,miter limit=1] (90-\bbeta:\T)coordinate(bbeta)--(0,0)--(90-\aalpha:\R)coordinate(aalpha);
		
		\draw[->] (90:\T-0.8)arc(90:90-\bbeta:\T-0.8);
		\coordinate (bet) at(90-0.5*\bbeta:\T-0.66){};
		\draw[->] (90:\R-0.5)arc(90:90-\aalpha:\R-0.5);
		\coordinate (alp) at(90-0.5*\aalpha:\R-0.35){};

\end{scope}
\end{scope}

%%%%% F R O N T %%%%%
\begin{scope}[yshift=\shiftX*-\xComplement,xshift=\shiftX]
\begin{scope}[every node/.append style={yslant=\xOffset,xscale= \xScale},yslant=\xOffset,xscale= \xScale,xshift=0]

		\foreach \m in {-3,...,4}{
			\foreach \n in {-2,...,1}{
			\draw[opacity=0.25,miter limit=1] (bbeta)--(-0.25+0.5*\m,0.25+\n*0.5*0.85)--(aalpha);
			}
		}
	
\end{scope}
\end{scope}
	
	\node at(alp){$\alpha$};
	\node at(bet){$\beta$};
	
	\draw[densely dashed] (bbeta)--coordinate[pos=0.82](temp)(aalpha);
	\draw (temp)--++(-10:0.35)+(0,-2pt)node[right,scale=0.8]{No LOS path (obstructed)};
	
	\node[circle,ball color=white,minimum size=11pt] at(bbeta){};
	\node[circle,ball color=white,minimum size=11pt] at(aalpha){};
	
	\draw (bbeta)+(0,-0.13)--++(0,0.1)--+(-0.1,0)--+(0,-0.133)--+(0.1,0)--+(0,0);
	\draw (aalpha)+(0,-0.13)--++(0,0.1)--+(-0.1,0)--+(0,-0.133)--+(0.1,0)--+(0,0);
	\node[scale=0.66,above=7pt] at(bbeta){Tx antenna};
	\node[scale=0.66,below=7pt] at(aalpha){Rx antenna};
	
	\path (aalpha)--coordinate[pos=0.05](tempa)coordinate[pos=0.1](tempb)(0,0);
	\draw[<-,line width=0.9pt] (tempa)--(tempb);
	
	\draw (0.75,1.4)--+(10:0.4)node[right,scale=0.9] {RIS (8$\times$4 elements)};

	\end{tikzpicture}
	\vspace*{-1.5mm}
	\caption{Geometric illustration of the setup: The RIS facilitates an NLOS communication channel between the transmitter and the receiver antennas, positioned at angles $\beta$ and $\alpha$ w.r.t. the surface normal, respectively.}%(at incident angle $\beta$) to the receiver (at the reflection angle $\alpha$).}%; the direct LOS path is assumed to be obstructed.}
	\label{fig:overall_setup}
	\vspace*{1mm}
\end{figure}

The wireless communication scenario under consideration is shown in \figref{fig:overall_setup}: A flat RIS consisting of several independently tunable antenna elements establishes a non-line-of-sight (NLOS) wireless channel from the transmitter antenna to the receiver antenna. The former is located at an incident angle $\beta$ w.r.t. the normal of the surface, whereas the latter is placed at the reflection angle $\alpha$. The line-of-sight (LOS) path between the transmitter and receiver is assumed to be obstructed, i.e., no LOS connection exists. This is not required, but a very practical scenario that makes the RIS more significant. 

%The distances of the transmitter and receiver could be different, but for simplicity and ease in practical verification, they are both always at a distance of $d=2$\,m from the center of the RIS. The overall methodology is not restricted to reflecting RIS, also transmitting RIS could be considered. However, the measurement setup restricts the considerations to incident and reflection angles in the range $\pm$90° off broadside from the RIS.

Generally, the transmitter and receiver antennas are located at different distances from the RIS. However, for simplicity and ease of practical verification, herein both are positioned at a constant distance $d=2$\,m from the RIS center. Moreover, the measurement setup confines considerations to incident and reflection angles within the $\pm$90° range off broadside from the RIS. %However, the methodology does not exclusively pertain to reflecting RISs; transmitting RISs or simultaneously transmitting and reflecting STAR-RIS~\cite{liu21star} could also be explored.
However, the methodology is not exclusive to reflecting RISs; it could also be applied to transmitting RISs or simultaneously transmitting and reflecting STAR-RISs~\cite{liu21star}.

%\subsection{RIS Design Aspects} % Considerations

\subsection{Review of Convex Optimization of RIS-based Systems}%{Review of RIS Optimization via Relaxation-Based Maximum Power Transfer Method}%Previously Proposed Optimization Method}
The methodology proposed in \cite{lang23} treats the entire wireless link from the transmitter antenna via RIS to the receiver antenna as an $(N+2)\times (N+2)$ impedance matrix $\bs{Z}$, where $N$ represents the number of antennas (ports) of the RIS that contain tunable elements. The unloaded impedance matrix of the entire system is obtained via simulation, e.g., in ANSYS HFSS with FEBI boundaries to cut down simulation time. The LOS path is omitted by zeroing the respective matrix entries.

Maximizing the power at the receiver antenna port, while holding the transmit power constant and enforcing passivity at all RIS antenna ports, yields optimal reactance values for these ports. Hence, when the RIS antenna ports are then terminated with these optimal reactances, the optimal operation conditions are attained and the RIS provides the maximum performance enhancement of the wireless channel in that configuration.

%As discussed in \cite{lang23}, this optimization relies on a convex semidefinite relaxation of the otherwise non-convex (i.e. "`unsolvable"') problem. Although there is no guarantee that this relaxation is tight, meaning necessarily leads to the true global optimum of the actual problem, there is a simple test that can confirm this. For all test cases, the optimization was successful.

As outlined in \cite{lang23}, this optimization method employs a convex semidefinite relaxation of the otherwise non-convex and thus "unsolvable" problem. Although there is no assurance that this relaxation is tight\,---\,meaning it necessarily leads to the true global optimum of the actual problem\,---\,a straightforward test can confirm its validity. As experienced in the original paper, the optimization was successful in all tested cases.

\vspace*{4pt}
\section{Practical Implementation}

\subsection{Reflectarray Design}
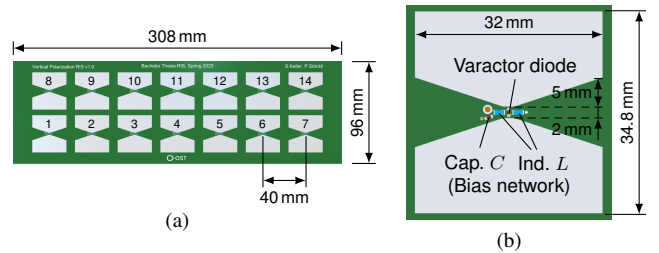
\begin{figure}[!htb]\centering
	\vspace*{-2mm}
	\begin{tikzpicture}[scale=0.9]
	
		\normalfont\footnotesize
		\node[yslant=0] at(-4.85,-1.6){(a)};
		\node[yslant=0] at(0,-1.9){(b)};
		
		\begin{scope}[yslant=0,every node/.append style={yslant=0},>=latex]\sffamily\small
		%\draw (0,0)rectangle(7,2);
		%\draw (0,1)--(7,1);
		%\draw (1,0)--(1,2) (2,0)--(2,2) (3,0)--(3,2) (4,0)--(4,2) (5,0)--(5,2) (6,0)--(6,2);
		\begin{scope}[xshift=-4.85cm]
		\node at(-0,0){\includegraphics[width=4.32cm]{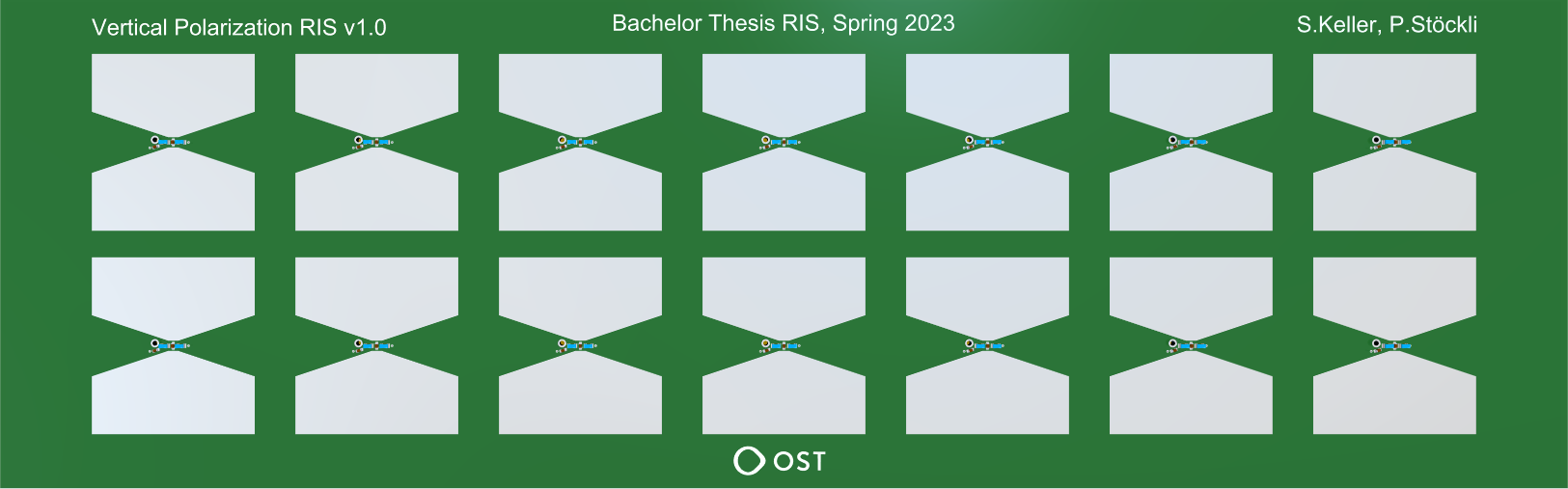}};
		\draw[<->] (-2.4,0.95)--node[above,scale=0.7]{308\,mm}(2.4,0.95);
		\draw (-2.4,1.05)--(-2.4,0.85) (2.4,1.05)--(2.4,0.85);
		\draw[<->] (2.8,-0.75)--node[above,scale=0.7,rotate=90]{96\,mm}(2.8,0.75);
		\draw (2.55,-0.75)--(2.9,-0.75) (2.55,0.75)--(2.9,0.75);
		\draw[<->] (1.25,-1)--node[below=2pt,scale=0.7]{40\,mm}(1.88,-1);
		\draw (1.25,-1.1)--+(0,0.75) (1.88,-1.1)--+(0,0.75);
		
		\foreach \n in {1,...,7}{
			\node[scale=0.5] at(\n*0.625-0.625-1.88,-0.17){\n};
		}
		\foreach \n in {8,...,14}{
			\node[scale=0.5] at(\n*0.625-0.625*8-1.88,0.47){\n};
		}
		
		\end{scope}
		
		\node at(0,0){\includegraphics[width=2.7cm]{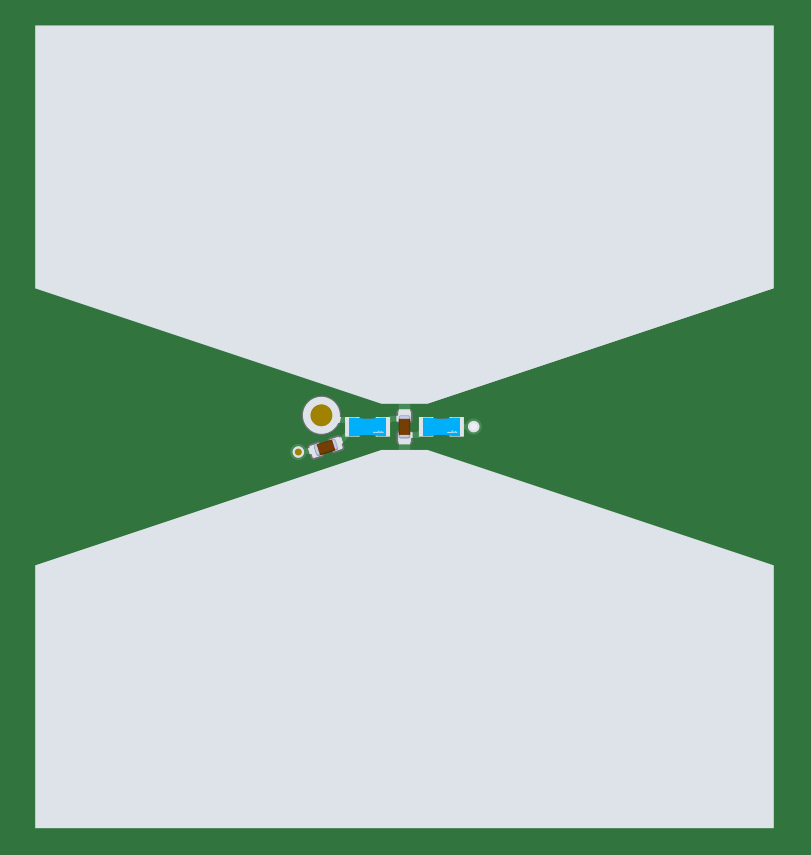}};
		
		\draw (0,0)--+(80:0.5)node[above,scale=0.8]{Varactor diode};
		\node[scale=0.8,align=left] at(0,-1.1){(Bias network)};
		\node[scale=0.8] (cap)at(-0.5,-0.75){Cap. $C$};
		\node[scale=0.8] (ind)at(0.5,-0.75){Ind. $L$};
		\draw (cap)--(-0.3,-0.08);
		\draw (ind)--(-0.15,-0.01);
		\draw (ind)--(0.15,-0.01);
		%\draw (-0.3,-0.12)--+(-100:0.8)node[below]{ cap. $C$};
		%\draw (-0.15,-0.01)--+(-70:0.4)node[below]{Bias ind. $L$};
		
		\draw[densely dashed] (0.075,0.08)--(1.4,0.08);
		\draw[densely dashed] (0.075,-0.08)--(1.4,-0.08);
		\draw[<->] (1.3,0.08)--node[left,scale=0.7]{5\,mm}(1.3,0.5);
		\draw (1.37,0.5)--(1.37-0.13,0.5);
		\draw[<-] (1.3,-0.08)--node[left,scale=0.7]{2\,mm}+(0,-0.33);
		\draw (1.3,-0.08)--(1.3,0.08);
		
		%\draw[<->] (-1.5,1.75)--node[above,scale=0.7]{32\,mm}(1.5,1.75);
		%\draw (-1.5,1.65)--(-1.5,1.75+0.1) (1.5,1.65)--(1.5,1.75+0.1);
		\draw[<->] (-1.37,1.15)--node[above,scale=0.7]{32\,mm}(1.37,1.15);
		\draw[<->] (1.9,1.5)--node[above,scale=0.7,rotate=90]{34.8\,mm}(1.9,-1.5);
		\draw (1.4,1.5)--(2,1.5) (1.4,-1.5)--(2,-1.5);
		
		\end{scope}
	\end{tikzpicture}
	\vspace*{-3mm}
	\caption{The full RIS (a) and its unit cell (b) with the varactor diode and bias network consisting of two inductors and a capacitor shown.}
	\label{fig:unit_cell}
	\vspace*{-3mm}
\end{figure}

The RIS under consideration follows the design from \cite{lang23}. However, for simplicity it only consists of two rows of 7 flat dipoles, i.e., a total of 14 tunable elements, as shown in \figref{fig:unit_cell}. The frequency under consideration is 3.55\,GHz, the center of the 5G band “n78”. The antenna elements are spaced 40\,mm apart from one another and 1.6\,mm-thick copper-backed FR4 is used as dielectric substrate. %\figref{fig:controller} shows the fully assembled RIS. %The fully assembled RIS can be seen in 

\subsection{Varactor Diodes}
Reverse-biased varactor diodes are used as tunable capacitors. While this is common for applications below optical and mmWave frequency ranges, also PIN diode designs can be found in recent literature, see for example \cite{gharbieh23,oh23}. Beyond that, MEMS and liquid crystals have also been investigated for tunability \cite{rana23}. Note that the method used herein is not limited to varactor diodes; any tunable element that can be operated to provide a particular reactance value is suitable.

\begin{figure}[!b]\centering
	\vspace*{-5.5mm}
	\begin{tikzpicture}\sffamily\small
		\node at(0,0){\includegraphics[width=8.2cm,clip,trim={5mm 0mm 5mm 0mm}]{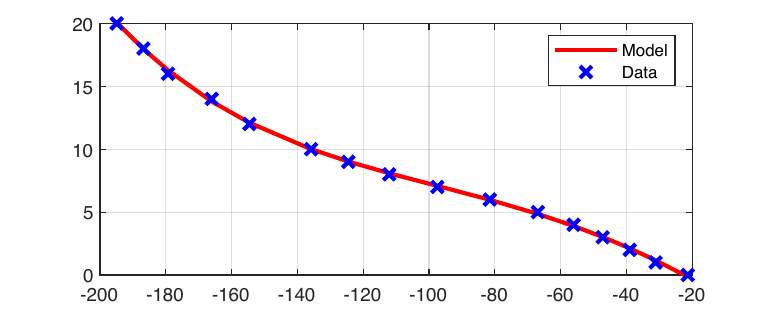}};
		\node[rotate=90,scale=0.9] at(-4,0){Voltage $v$ in V};
		\node[scale=0.9] at(0,-2.2){Reactance $x$ (capacitive) in $\Omega$};
	\end{tikzpicture}
	\vspace*{-2mm}
	\caption{The Skyworks SMV2201-040LF varactor diode.}%: values from the datasheet and the model.}
	\label{fig:diode_model}
	%\vspace*{-2.5mm}
\end{figure}

The SMV2201-040LF varactor diode from Skyworks \cite{SMV2201} was selected for this project due to its capability to achieve one of the smallest capacitances available, ranging from 0.23 to 2.1\,pF as shown in  \figref{fig:diode_model},
 covering most of the desired values. The datasheet indicates that different devices can vary by 50\% for very small capacitances and 10\% for larger ones, which may impact the accuracy of the validation in practice.

\subsubsection{Varactor Diode Model}
In general, for reliable results, the equivalent RLC circuit of the varactor diodes must be extracted from measurements \cite{ratajczak23}. However, in this case, the limited values available in the datasheet were used to find a low-order polynomial for the required voltage $v$ to obtain a desired reactance value $x$. As shown in \figref{fig:diode_model}, the following polynomial approximation is used for the SMV2201-040LF:
%Simple polynomial approximation: The required voltage $v$ to obtain the desired reactance $X$ is obtained with sufficient accuracy as
\begin{equation}
	v(x) = a + bx + cx^2 +dx^3 + ex^4
\end{equation}
with the coefficients $a = -3.55$\,V, $b = -1.77\times10^{-1}\,\frac{\text{V}}{\Omega}$, $c=-8.28\times10^{-4}\,\frac{\text{V}}{\Omega^2}$, $d = 8.5\times 10^{-8}\,\frac{\text{V}}{\Omega^3}$, and $e = 1.47\times 10^{-8}\,\frac{\text{V}}{\Omega^4}$.

The loss resistance of the diode is about $5.4\,\Omega$ according to the datasheet \cite{SMV2201}. This was confirmed by measurement and cannot be neglected when using it for large capacitance values.

\subsubsection{Bias Network}
%The bias network for the varactor diodes consists of two inductors with $L=33$\,nH $\pm 2$\% (Murata LQW18AN33NG00D) and a capacitor of 3.9\,pF $\pm1$\,pF (Johanson xxyy). They were both chosen so that their self resonance frequencies lie close to the operating frequency. As can be seen in \figref{fig:unit_cell}, the layout was designed to maintain symmetry as much as possible.

The bias network, designed for the varactor diodes, incorporates two inductors of $L=33$\,nH $\pm 2$\% (Murata LQW18AN33NG00D) and a capacitor, with a capacitance of $C=3.9$\,pF $\pm0.1$\,pF (Johanson 500R07S39RBV4S). Both components were strategically selected to ensure that their self-resonance frequencies align closely with the operational frequency. As shown in \figref{fig:unit_cell}(b), the layout was designed to maintain symmetry as much as possible.
%The layout, meticulously designed to preserve as much symmetry as possible, is illustrated in \figref{fig:unit_cell},.

\subsubsection{Limited Range of Values}
The optimization process may potentially yield either inductive or capacitive reactance values that are unattainable with the varactor diodes. In such instances, the corresponding ports are terminated with the nearest realizable suboptimal values\,---\,specifically, approximately 2.1\,pF (equivalent to a reactance $x\approx-20\,\Omega$) at 0\,V or 0.23\,pF ($x\approx -200\,\Omega$) at 20\,V. The remaining impedance matrix could be re-optimized utilizing the previously mentioned optimization framework or the optimization could be adapted to constrain the reactive loads. However, this is beyond the scope of this paper.% This approach ensures the attainment of the best performance from the RIS that is realizable with the selected varactor diodes.

\subsection{RIS Controller}
\begin{figure}[!b]\centering\sffamily\small\vspace*{-2.5mm}
	\begin{tikzpicture}
		\node at(0,0){\includegraphics[width=7.5cm,clip,trim={0 1mm 0 7mm}]{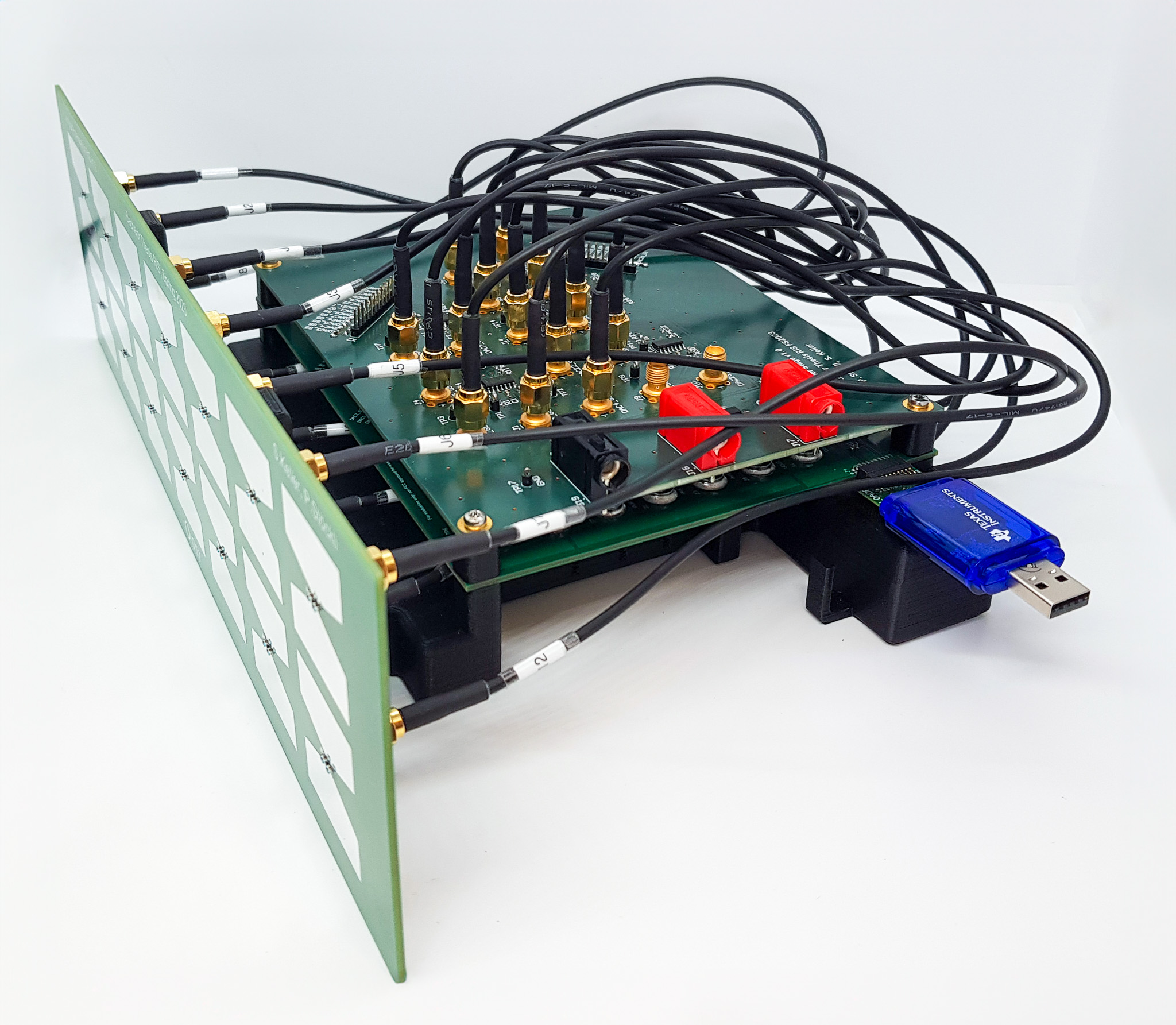}};
		\draw (-2.1,-1.05)--+(-120:1)node[text width=2.5cm,align=center,scale=0.9,below]{RIS\\ (7$\times$2 elements)};
		\draw (0,-0.15)--+(-80:1.4)node[below,text width=2cm,align=center,scale=0.9]{DAC60096 evaluation board};
		\draw (-1.5,1.7)--+(100:0.68)node[above,text width=2cm,align=center,scale=0.9]{Amplifier board};
		\draw (2.8,-0.45)--+(-90:0.6)node[below,text width=2cm,align=center,scale=0.9]{USB to SPI adapter};
	\end{tikzpicture}
	\vspace*{-2.5mm}
	\caption{Photo of the fully assembled RIS featuring the controller, composed of the DAC and amplifier boards, mounted on its rear side.}
	\label{fig:controller}
	%\vspace*{-2.5mm}
\end{figure}

An RIS controller board was developed to supply the required voltages to the varactor diodes. Anticipating the needs of larger designs in future projects, the DAC60096 by Texas Instruments was chosen. Because its output voltage is limited to 10.5\,V, but 20\,V are required to obtain sufficiently low capacitance values from the varactor diodes, it is utilized in conjunction with an amplifier stage, depicted in \figref{fig:controller}. 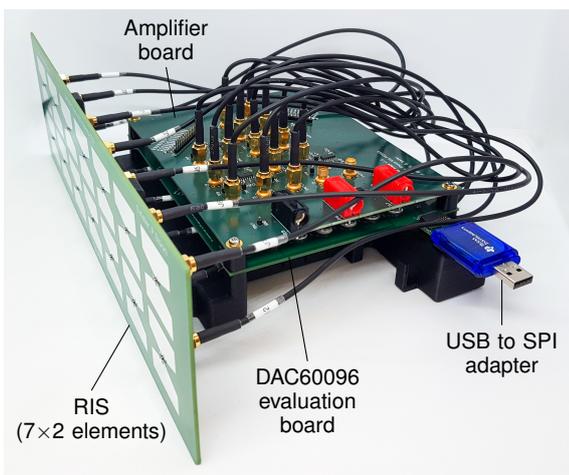

The NCS21914 operational amplifiers by ON Semiconductors were selected for this stage, due to their ability to deliver satisfactory gain and stability while preserving low-power consumption. Additionally, they provide a high input impedance and exhibit rail-to-rail capabilities when deployed with an asymmetric power supply.

%To provide the necessary voltages to the varactor diodes, an RIS controller board was developed. To be ready for larger designs in the future, the DAC60096 by Texas Instruments was selected. Since its output voltage is limited to 10.5\,V but 20\,V are required to achieve sufficiently low capacitance values from the varactor diodes, it is used in conjunction with an amplifier stage, shown in \figref{fig:controller}. NCS21914 operational amplifiers by ON Semiconductors were selected for this stage due to their capacity to provide adequate gain and stability while maintaining low-power consumption, ensuring a high input impedance and its rail-to-rail capabilities when using an asymmetric power supply.

% is used in conjunction with an amplifier stage, shown in \figref{fig:controller}. The addition of an amplifier stage is vital because the DAC is limited to outputting voltages of up to 10.5\,V, while a minimum of 20\,V is required to achieve sufficiently low capacitance values from the varactor diodes. NCS21914 operational amplifiers by ON Semiconductors were selected for this stage due to their capacity to provide adequate gain and stability while maintaining low-power consumption and ensuring a high input impedance.

\vspace*{4pt}
\section{Measurements}

\subsection{Measurement Setup}
%Measuring RISs is not a trivial task, since they often require distances that go beyond small anechoic chambers. This was the case here as well. Therefore, three movable absorber walls were constructed and a measurement setup consisting of two movable arms and a central pivot point was developed. A Raspberry Pi controls the two 24\,V motors that move the arms independently until the desired angles are reached. By mounting the antennas on these arms, the antennas are always aimed at the center of the RIS. It is illustrated in \figref{fig:measurement_setup}(a) and a photo is shown in (b). 

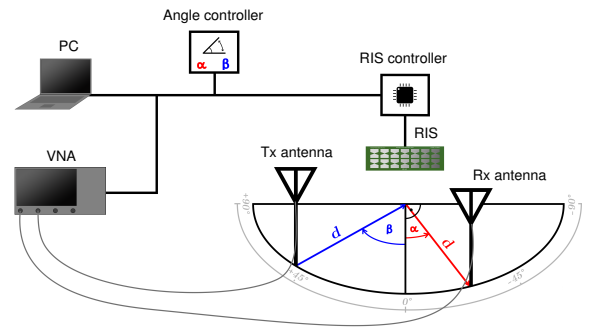
\begin{figure}[!htb]\centering
	\vspace*{-1mm}
	\begin{tikzpicture}\sffamily\small
		\node at(0,0){\includegraphics[width=7.5cm,clip,trim={5mm 0 2mm 10mm}]{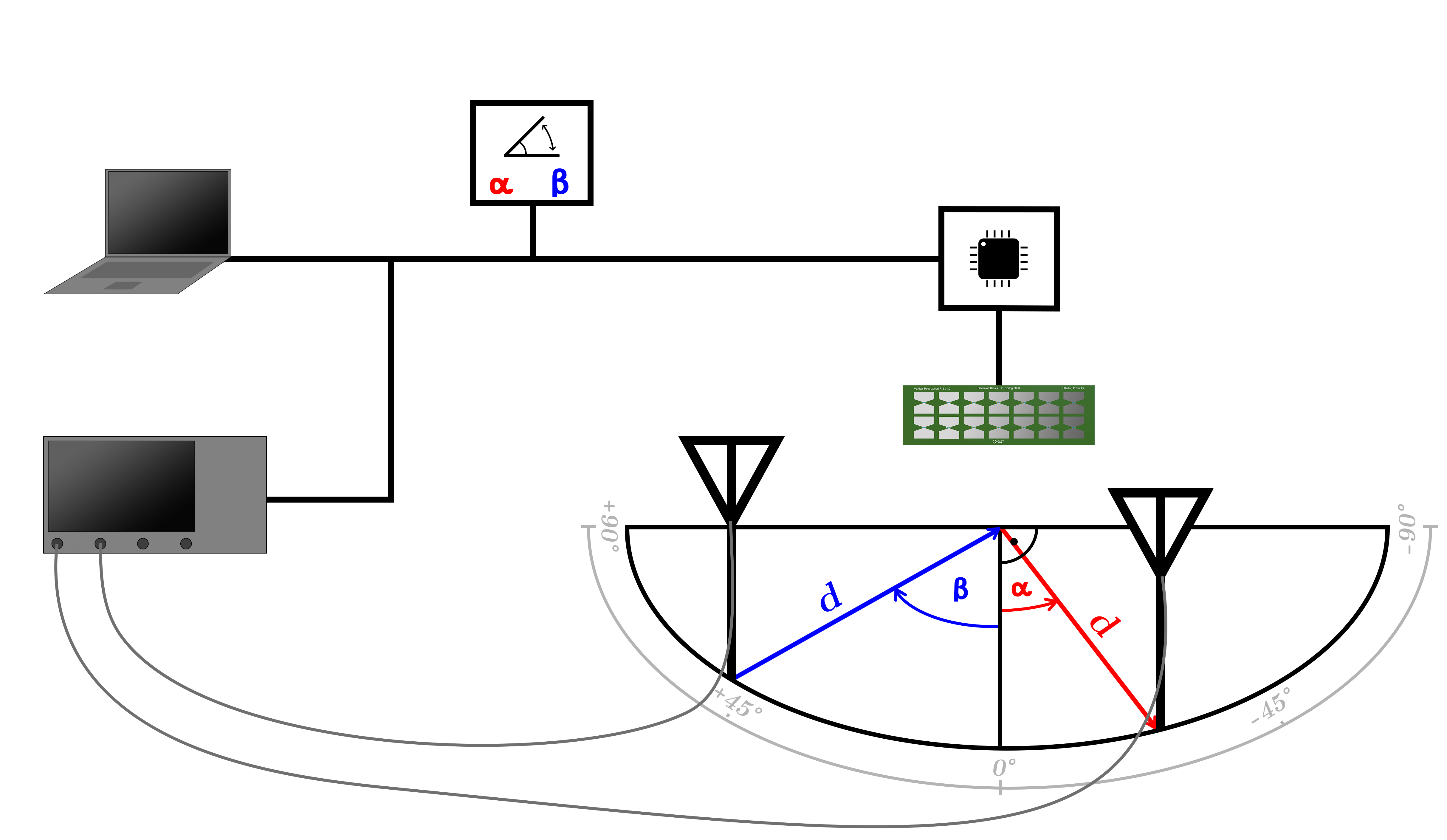}};
		\node[scale=0.6] at(-1.1,2.1){Angle controller};
		\node[scale=0.6] at(1.4,1.55){RIS controller};
		\node[scale=0.6] at(-3,1.7){PC};
		\node[scale=0.6] at(-3.1,0.25){VNA};
		\node[scale=0.6] at(0.0,0.3){Tx antenna};
		\node[scale=0.6] at(2.8,0.0){Rx antenna};
		\node[scale=0.6] at(1.7,0.55){RIS};
	\end{tikzpicture}	
	%\subfloat[]{\begin{tikzpicture}
		%\draw (0,0)rectangle(8,5);
		%\draw (0,0)--(8,5);
	%\end{tikzpicture}}
	\vspace*{-4mm}
	\caption{The RIS measurement setup: The two movable arms and the RIS are operated via a PC, which also collates measurement data from the VNA.}
	%\caption{RIS measurement setup: PC-operated arms and controller, and data acquisition from the VNA.}
	\label{fig:description_measurement_setup}
	\vspace*{0.5mm}
\end{figure}

Measuring RISs is inherently challenging, because it often requires distances that exceed the capabilities of smaller anechoic chambers\,---\,a scenario encountered in this study as well. To address this, three movable absorber walls were constructed and a dedicated measurement setup was developed, as illustrated in  \figref{fig:description_measurement_setup}. It features two adjustable arms anchored on a central pivot point where their angles are measured by Amphenol Piher angle sensors (PST-360GB1A-C0000-ERA-05K). The arms are independently moved to desired angles by two 24\,V DC motors, all controlled by a Raspberry Pi 4 (“Angle controller”) connected to a PC over Ethernet.
%\begin{figure}[!htb]\centering
%	\vspace*{-3mm}
%	\subfloat[]{\includegraphics[width=7.75cm]{measurement_setup.png}}\\[3mm]
%	\subfloat[]{\includegraphics[width=8.25cm,clip,trim={0 0 2mm 3mm}]{DSC08677-Enhanced-NR.jpg}}
%	%\subfloat[]{\begin{tikzpicture}
%		%\draw (0,0)rectangle(8,5);
%		%\draw (0,0)--(8,5);
%	%\end{tikzpicture}}
%	\caption{Description (a) and photo (b) of the entire measurement setup.}
%	\label{fig:measurement_setup}
%\end{figure}

Since the received signal of interest is of very low power, reflections from the environment can notably distort the measurement results through interference. Therefore, two QRH11 quad ridged horn antennas by RFspin are used as transmitter and receiver antennas. They provide about 11\,dBi gain at the frequency of interest and allow for all polarizations to be measured simultaneously, which is not required here, but useful for future considerations. These antennas are mounted on the arms and consistently point to the center of the RIS.
\begin{figure}[!htb]\centering
	%\vspace*{-3mm}
	\includegraphics[width=8.25cm,clip,trim={0 2mm 4mm 1mm}]{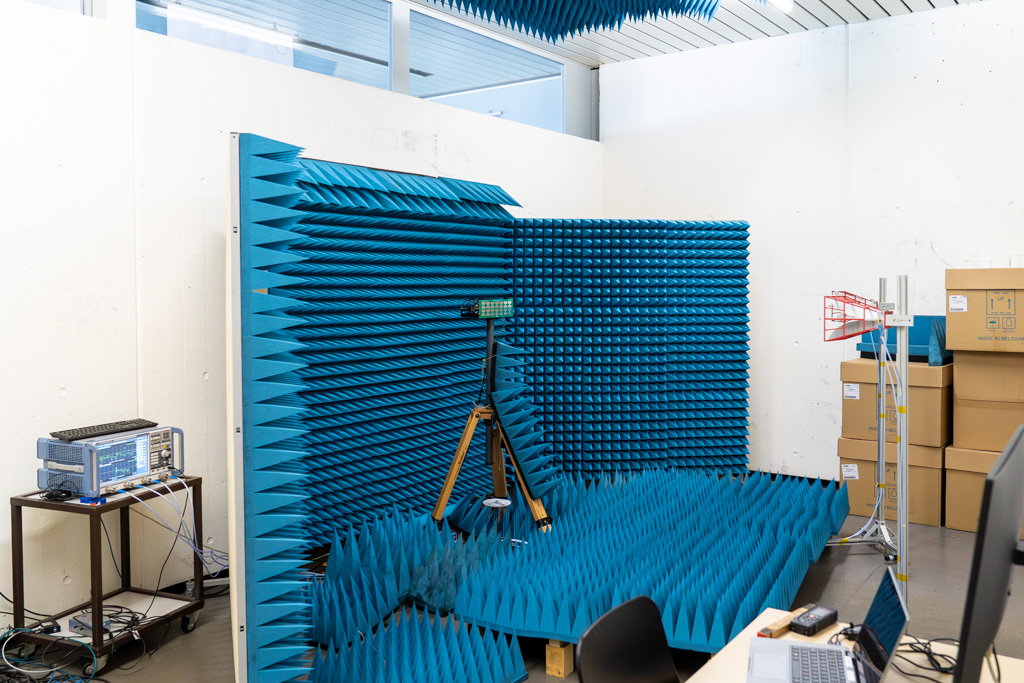}
	%\subfloat[]{\begin{tikzpicture}
		%\draw (0,0)rectangle(8,5);
		%\draw (0,0)--(8,5);
	%\end{tikzpicture}}
	\vspace*{-2mm}
	\caption{Photo of the measurement setup, featuring two movable antenna arms, a VNA, and the mounted RIS, surrounded by numerous absorbers.}
	%\caption{Photo of the measurement setup consisting of three movable absorber walls, two selectronically controllable antenna arms, a VNA, and the mounted RIS.}
	\label{fig:photo_measurement_setup}
	\vspace*{-2mm}
\end{figure}%Moreover, this manufacturer provides measured gains over the entire spectrum for each antenna, which is quite helpful. %Because in this case these antennas are always aimed at the RIS and the direct path between the antennas is omitted, this value is sufficient.

This setup has demonstrated significant utility and versatility. However, it presents the following notable limitations:
\begin{itemize}
\item \textit{Reflections from the Environment:}\enspace 
In the measurement setup shown in \figref{fig:photo_measurement_setup}, effectively eliminating all environmental reflections proves to be a significant challenge. %The primary issue stems from reflections arriving close to the time gate, especially those paths near a distance of $2d$, equivalent to the distance of the signal of interest. 
Additional absorbers are placed in the room, particularly on the floor and on the ceiling, to mitigate reflections experienced from those locations. The use of time gating, described in the next subsection, helps further.% particular distance. 

\item \textit{Far-Field Condition:}\enspace 
At a distance of 2\,m, the transmitter and receiver antennas are situated approximately at the threshold of the RIS's far field. Moreover, at small angles $\abs{\alpha+\beta}\leq15$°, the two quad ridge antennas are within their respective near fields. However, such small angles are not possible with the setup anyway. These conditions may impact the measurements, but since the direct path is omitted, the influence is expected to be negligible.

\item \textit{Measurement Angle Uncertainty:}\enspace 
The angle controller within the measurement setup exhibits an absolute angle error of $\pm$1°, introducing a potentially non-negligible factor of uncertainty in the reproducibility of measurement results. Moreover, even minor deviations in the antenna angle can yield disparate measurement outcomes. %To facilitate a comparative analysis of the measurement results under varied settings (e.g. loaded/unloaded), successive measurements are conducted at specific angles. %A notable challenge arises as the performance of the capacitances, obtained from measurements in Section 6.1, cannot be precisely evaluated at their optimized points due to this inherent inaccuracy. Nevertheless, this is anticipated to induce only a marginal distortion along the x-axis during a sweep across a comprehensive range of angles.
\end{itemize}

\subsection{Time Gating}
%Time gating is a well-known technique for VNA measurements often employed to mitigate imperfections due to reflections of the environment \cite{gonzalez16,phumvijit17,williams90}. 

%The gate was applied so as to be 1\,ns before the expected reflection and 10\,ns wide. Hence, it is used to omit the direct LOS path, and reduce inaccuracies due to other reflections of the environment. 

Time gating is a well-established technique in Vector Network Analyzer (VNA) measurements. For antenna measurements, it is frequently utilized to mitigate imperfections arising from environmental reflections \cite{williams90,gonzalez16,phumvijit17}, by applying a filter that basically omits signals outside a certain time-frame. %Time gating is a well-established technique in Vector Network Analyzer (VNA) measurements. For antenna measurements, it is frequently utilized to mitigate imperfections .

\begin{figure}[!htb]\centering
	\hspace*{-0.25mm}\includegraphics[width=1\columnwidth,clip,trim={20.1mm 24mm 20mm 24mm}]{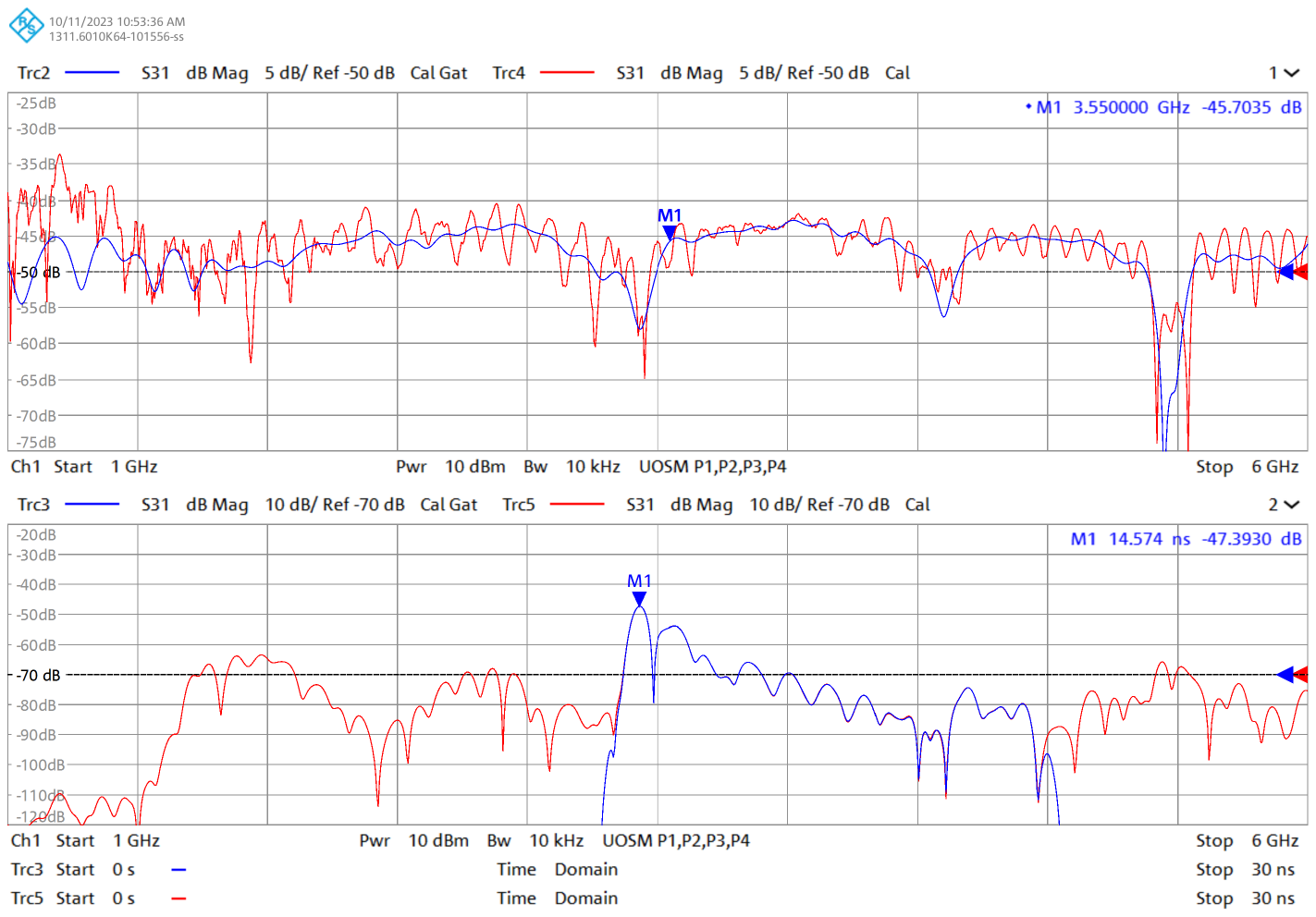}
	%\begin{tikzpicture}
	%	\draw (0,0)rectangle(8,5);
	%	\draw (0,0)--(8,5);
	%\end{tikzpicture}
	\vspace*{-7mm}
	\caption{Example of time gating (unmodified screenshot) in both frequency domain (top) and time domain (bottom): The time gate starting at 13\,ns omits the direct (LOS) path arriving around 5\,ns, while accounting for the indirect path via RIS peaking around 14\,ns.  This case corresponds to $\alpha=20$° and $\beta=-10$° with the original shown in red and the time-gated in blue.}
	\label{fig:timegating}
	\vspace*{-2mm}
\end{figure}

In this study, a Hann-type gate (i.e., filtering window) was set to start approximately 1\,ns prior to the anticipated reflection off the RIS, with a width of 10\,ns. Consequently, this effectively omits the direct line-of-sight (LOS) path and minimizes inaccuracies attributed to other environmental reflections. \figref{fig:timegating} shows a typical result. Note that the gate/window width should not be chosen too small; the reflection off the RIS takes considerably longer than from a plane reflector.

\subsection{Verification Results and Discussion}

\subsubsection{Maximum Bistatic Radar Cross Section (BRCS)}

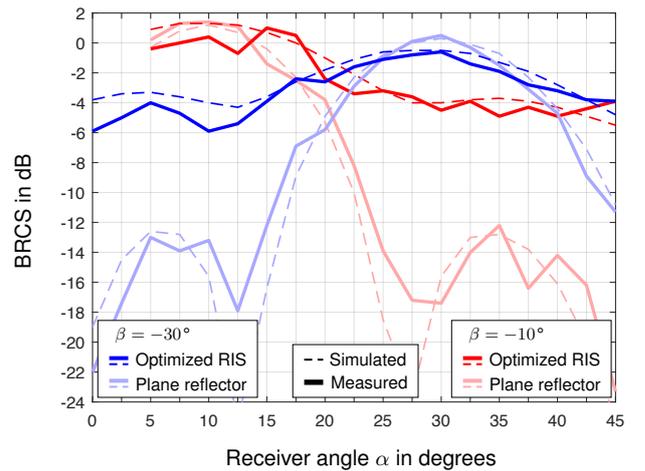
\begin{figure}[!b]\centering
	\vspace*{-4mm}
	\begin{tikzpicture}\sffamily\small
		\node at(0,0){\includegraphics[width=8.0cm,clip,trim={7mm 5mm 7mm 3mm}]{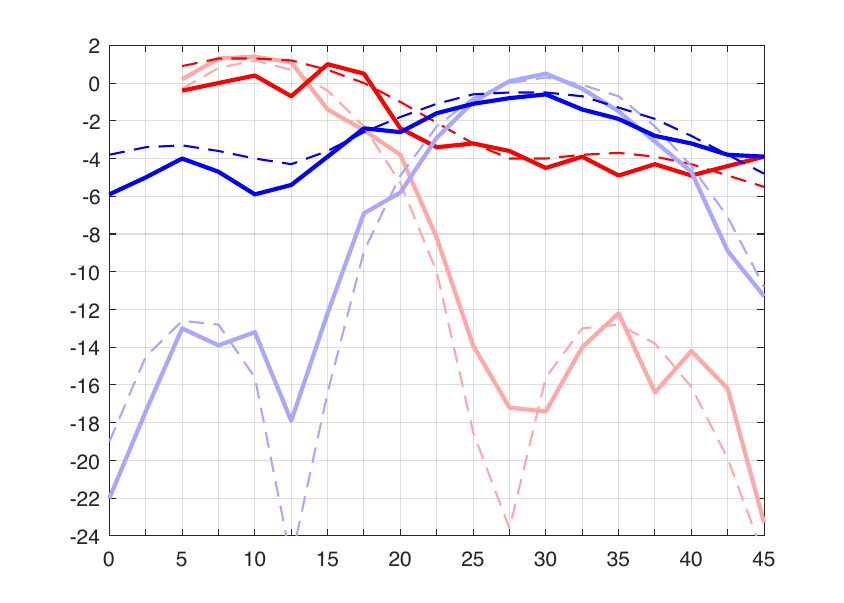}};
		\node[scale=0.9] at(0.25,-3.3){Receiver angle $\alpha$ in degrees};
		\node[scale=0.9,rotate=90] at(-4.2,0){BRCS in dB};
		
		%\begin{scope}[xshift=-3cm,yshift=-0.7cm]
			%\def\d{0.3}
			%\def\tscale{0.72}
			%\draw[fill=white] (-0.15,0.2)rectangle(3.64,-5*\d-0.2);
			%\draw[line width=1pt,red,densely dashed] (0,0)--+(0.25,0)node[right,text=black,scale=\tscale]{Optimum simulated $\beta=-10$°};
			%\draw[line width=2pt,red] (0,-\d)--+(0.25,0)node[right,text=black,scale=\tscale]{Optimum measured $\beta=-10$°};
			%\draw[line width=1pt,blue,densely dashed] (0,-2*\d)--+(0.25,0)node[right,text=black,scale=\tscale]{Optimum simulated $\beta=-30$°};
			%\draw[line width=2pt,blue] (0,-3*\d)--+(0.25,0)node[right,text=black,scale=\tscale]{Optimum measured $\beta=-30$°};
			%\draw[line width=1pt,black,densely dashed] (0,-4*\d)--+(0.25,0)node[right,text=black,scale=\tscale]{Plane reflector sim. $\beta=-10$°};
			%\draw[line width=2pt,black] (0,-5*\d)--+(0.25,0)node[right,text=black,scale=\tscale]{Plane reflector meas. $\beta=-10$°};
		%\end{scope}
		
			\def\tscale{0.72}
			\def\d{0.3}
			
			\begin{scope}[xshift=-3.065cm,yshift=-1.65cm,line width=0.25pt]
				\draw[fill=white] (-0.15,0.2)rectangle(1.95,-2*\d-0.22);
				\node[right,text=black,scale=\tscale] at(0,0){$\beta=-30$°};
				\draw[line width=1.5pt,blue] (0,-\d)--++(0.25,0)+(0,-0.04)node[right,text=black,scale=\tscale]{Optimized RIS};
				\draw[line width=0.75pt,blue,densely dashed] (0,-\d-0.08)--++(0.25,0);
				
				\draw[line width=1.5pt,blue!33!white] (0,-2*\d)--++(0.25,0)+(0,-0.04)node[right,text=black,scale=\tscale]{Plane reflector};
				\draw[line width=0.75pt,blue!33!white,densely dashed] (0,-2*\d-0.08)--++(0.25,0);
			\end{scope}
			
			\begin{scope}[xshift=1.59cm,yshift=-1.65cm,line width=0.25pt]
				\draw[fill=white] (-0.15,0.2)rectangle(1.95,-2*\d-0.22);
				\node[right,text=black,scale=\tscale] at(0,0){$\beta=-10$°};
				\draw[line width=1.5pt,red] (0,-\d)--++(0.25,0)+(0,-0.04)node[right,text=black,scale=\tscale]{Optimized RIS};
				\draw[line width=0.75pt,red,densely dashed] (0,-\d-0.08)--++(0.25,0);
		
				\draw[line width=1.5pt,red!33!white] (0,-2*\d)--++(0.25,0)+(0,-0.04)node[right,text=black,scale=\tscale]{Plane reflector};
				\draw[line width=0.75pt,red!33!white,densely dashed] (0,-2*\d-0.08)--++(0.25,0);
			\end{scope}
		
			\begin{scope}[xshift=1cm,yshift=-1.97cm,line width=0.25pt]
			\def\d{0.3}
			\draw[fill=white] (-1.5-.15,0.2)rectangle(0,-1*\d-0.2);
			%\draw (0,0.2)--coordinate(temp)(0,-1*\d-0.2);
			%\draw (temp)--+(5,0);
			%\node[right,text=black,scale=\tscale] at(0,0){$\beta=-10$°};
			%\draw[line width=2pt,red] (1.3,0)--+(0.25,0)node[right,text=black,scale=\tscale]{Optimized RIS};
			%\draw[line width=2pt,red!33!white] (3.5,0)--+(0.25,0)node[right,text=black,scale=\tscale]{Plane refl. };
			%\draw[line width=2pt,blue!33!white] (3.5,-\d)--+(0.25,0)node[right,text=black,scale=\tscale]{Plane refl.};
			\draw[line width=0.75pt,black,densely dashed] (-1.5,0)--+(0.25,0)node[right,text=black,scale=\tscale]{Simulated};
			\draw[line width=2pt,black] (-1.5,-\d)--+(0.25,0)node[right,text=black,scale=\tscale]{Measured};
		\end{scope}
	\end{tikzpicture}
	\vspace*{-2mm}
	\caption{Resulting bistatic radar cross sections (BRCSs) of the optimized RIS and a plane reflector (lighter colors) for the same transmitter angle $\beta=-10$° (red) or $-30$° (blue) and various receiver angles $\alpha=0\,...\,45$° off broadside (solid: measured, dashed: simulated).}
	\label{fig:results}
	%\vspace*{-1mm}
\end{figure}

Similar to \cite{lang23}, the bistatic radar cross section (BRCS) is used as a figure of merit to validate the optimization procedure. \figref{fig:results} depicts the results for a transmitter incident angles $\beta=-10$° (red) and $\beta=-30$° (blue), respectively, across various angles $\alpha$ ranging from 0° or 5° to 45° at the aforementioned distance of 2\,m. Notably, there is good agreement between the simulated (dashed) and measured (solid lines) values, despite the challenges and influential factors discussed previously. Additionally, while the measurement setup proved to be quite sturdy, it is observed that even minor misalignment can cause deviations of $\pm1$\,dB in sensitive areas.%As can be seen, the agreement is quite good, despite all the challenges mentioned.

The lighter lines in \figref{fig:results} provide a comparison of the simulated and measured BRCS of a plane copper reflector with the same dimensions as the array for the two transmitter angles. Given that the simulation utilizes plane-wave excitation, nulls are anticipated to be less pronounced in the measurements.  The overall agreement is quite good, apart from minor discrepancies around 37.5°. %Observe the distinction between the dashed and solid lines: It is noteworthy that even such a small RIS can exert a substantial impact at angles sufficiently far from the specular direction.
It is worth noting the distinction between the red and black lines: Even such a small RIS can create a substantial impact at angles sufficiently far from the specular direction.

\subsubsection{Optimality of the Obtained Reactances}
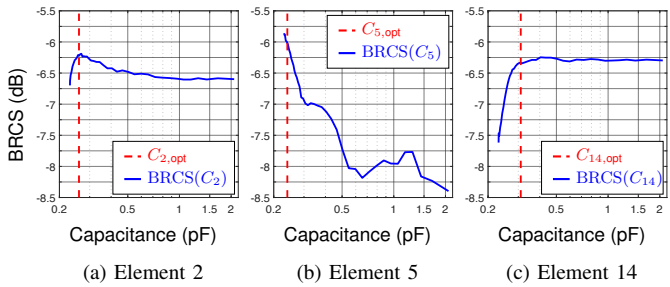
\begin{figure}[!htb]\centering
	\begin{tikzpicture}
		\def\d{2.82}
		\node[inner sep=0pt] at(0,0){\includegraphics[width=2.72cm,clip,trim={1mm 0.35mm 6.3mm 4mm}]{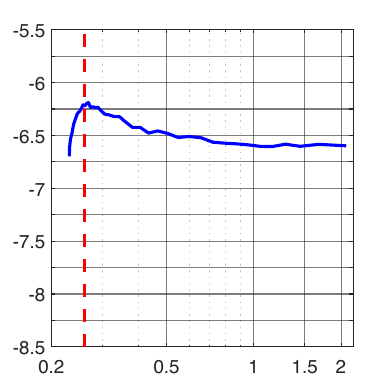}};
		\node[inner sep=0pt] at(\d,0){\includegraphics[width=2.72cm,clip,trim={1mm 0.35mm 6.3mm 4mm}]{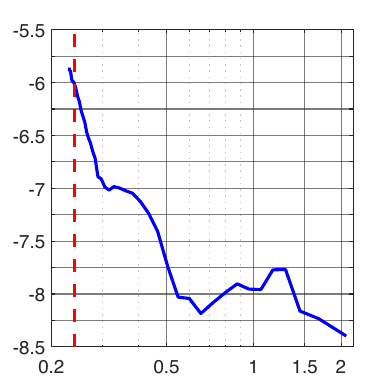}};
		\node[inner sep=0pt] at(2*\d,0){\includegraphics[width=2.72cm,clip,trim={1mm 0.35mm 6.3mm 4mm}]{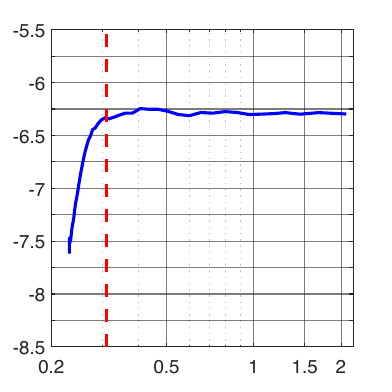}};
		\footnotesize
		\node at(0.1,-2.1){(a) Element 2};
		\node at(0.1+\d,-2.1){(b) Element 5};
		\node at(0.1+2*\d,-2.1){(c) Element 14};
		\sffamily\small
		\node[scale=0.8] at(0.1,-1.6){Capacitance (pF)};
		\node[scale=0.8] at(0.1+\d,-1.6){Capacitance (pF)};
		\node[scale=0.8] at(0.1+2*\d,-1.6){Capacitance (pF)};
		\node[scale=0.8,rotate=90,inner sep=0pt] at(-1.55,0){BRCS (dB)};
		
		\begin{scope}[xshift=-1.5mm,yshift=-0.5mm]
		\sffamily\small
		\draw[fill=white,line width=0.25pt] (-0.15,0.2-0.5)rectangle(1.45,-0.5-0.5);
		\draw[red,densely dashed,line width=0.8pt] (0,-0.5)--+(0.22,0)node[right,scale=0.7]{$C_{2,\text{opt}}$};
		\draw[blue,line width=0.8pt] (0,-0.5-0.3)--+(0.22,0)node[right,scale=0.7]{$\operatorname{BRCS}(C_2)$};
		\end{scope}
		
		\begin{scope}[xshift=-1.5mm+\d cm,yshift=16.2mm]
		\sffamily\small
		\draw[fill=white,line width=0.25pt] (-0.15,0.2-0.5)rectangle(1.45,-0.5-0.5);
		\draw[red,densely dashed,line width=0.8pt] (0,-0.5)--+(0.22,0)node[right,scale=0.7]{$C_{5,\text{opt}}$};
		\draw[blue,line width=0.8pt] (0,-0.5-0.3)--+(0.22,0)node[right,scale=0.7]{$\operatorname{BRCS}(C_5)$};
		\end{scope}
		
		\begin{scope}[xshift=-1.5mm+2*\d cm,yshift=-0.5mm]
		\sffamily\small
		\draw[fill=white,line width=0.25pt] (-0.15,0.2-0.5)rectangle(1.45,-0.5-0.5);
		\draw[red,densely dashed,line width=0.8pt] (0,-0.5)--+(0.22,0)node[right,scale=0.7]{$C_{14,\text{opt}}$};
		\draw[blue,line width=0.8pt] (0,-0.5-0.3)--+(0.22,0)node[right,scale=0.7]{$\operatorname{BRCS}(C_{14})$};
		\end{scope}
	\end{tikzpicture}
	%\subfloat[]{\begin{tikzpicture}[scale=0.9]
		%\node[inner sep=0pt] at(0,0){\includegraphics[width=2.75cm,clip,trim={1mm 0.35mm 6.3mm 4mm}]{sensitivity_element2_log_3.pdf}};
	%\end{tikzpicture}}
	%\hfill
	%\subfloat[]{\begin{tikzpicture}[scale=0.9]
		%\node[inner sep=0pt] at(0,0){\includegraphics[width=2.75cm,clip,trim={1mm 0.35mm 6.3mm 4mm}]{sensitivity_element5_log_3.pdf}};
	%\end{tikzpicture}}
	%\hfill
	%\subfloat[]{\begin{tikzpicture}[scale=0.9]
		%\node[inner sep=0pt] at(0,0){\includegraphics[width=2.75cm,clip,trim={1mm 0.35mm 6.3mm 4mm}]{sensitivity_element14_log_3.pdf}};
	%\end{tikzpicture}}
	\vspace*{-3mm}
	\caption{BRCS as a function of different capacitances $C_n$, to illustrate typical optimalities and sensitivities of the results obtained from the convex optimization algorithm ($C_{n,\text{opt}}$, red dashed line), for $\beta=-10$° and $\alpha=45$°.}
	\label{fig:ris_first_measurements}
	\vspace*{-1mm}
\end{figure}

 \figref{fig:ris_first_measurements} shows three typical behaviors of the capacitor values, obtained by varying one capacitor while maintaining the others at their optimized value. Evidently, the values derived from the optimization algorithm (red dashed lines) are indeed optimal (a) or at least close to optimal (b) and (c). Additionally, it can be observed that the optimum is insensitive to certain capacitances, while exhibiting significant sensitivity to others.
 
\vspace*{4pt}
\section{Summary \& Conclusions}
Reconfigurable Intelligent Surfaces (RISs), expected to play an important role in future wireless systems, can be designed through various methodologies. Among them is a powerful and versatile optimization framework, treating the RIS as a reflectarray and determining the optimum reactive loads in the antenna elements to maximize the power transferred between transmitter and receiver antennas.%which approaches them as reflectarrays and finds the optimum reactive loads in the antenna elements to maximize the power transferred from a transmitter a receiver antenna.

This paper presents a small RIS design to verify a previously proposed optimization method to find the optimum reactive loads in measurements. This task is challenging for various reasons, including the spatial requirements of the measurement setup as well as the limitations and tolerances of the varactor diodes used as variable capacitors. Despite these challenges, results validate the optimization; %show that the optimization result is meaningful and realizable. 
the obtained capacitance values are optimal or close to it in practice. The use of time gating has proved pivotal, effectively reducing reflections and omitting the line-of-sight path.

%Finally, it is shown that the realized RIS works as intended and its effectiveness is verified. It is remarkable that even a small RIS can have a considerable effect at angles far from the specular reflection. Furthermore, the convenient adjustability of the loading reactances allows for precise fine-tuning, helping to bridge any gaps between simulation and real-world outcomes.

%In conclusion, the implemented RIS not only functions as anticipated but also underscores its efficacy through verification. Notably, even a diminutive RIS can exert a substantial impact at angles divergent from the specular reflection. Moreover, the facile adjustability of the loading reactances enables meticulous fine-tuning, adeptly bridging disparities between simulation and actual results.

In conclusion, the realized RIS not only operates as anticipated, but also serves as a testament to the validity and efficacy of the convex optimization method. It is remarkable that even a small RIS can have a considerable effect at angles far from the specular reflection. Furthermore, the convenient adjustability of the loading reactances allows for precise fine-tuning, helping to bridge any gaps between simulation and real-world outcomes. 

%Furthermore, the adaptable adjustability of the loading reactances facilitates precise fine-tuning, effectively bridging the gap between simulation and real-world results. This attests to the practicability of the optimization method, opening avenues for further exploration and application in future RIS designs and implementations.

\vspace*{4pt}
\section*{Acknowledgment}

This work was carried out as part of the collaborative CHIST-ERA project \textit{"Towards Sustainable ICT: Sparse Ubiquitous Networks based on Reconfigurable Intelligent SurfacEs (SUNRISE)"}. The authors of OST are supported by the Swiss National Science Foundation (SNF) with Grant 203784.

\vspace*{4pt}
\bibliographystyle{myIEEEtran}
\bibliography{EuCAP2024_RIS_paper_noapprox_revision.bbl}

%\balance

\end{document}